 
\input amstex
\input epsf
\loadmsbm
\documentstyle{gsm}
\NoBlackBoxes
\font\tit=cmr10 scaled\magstep3
\font\sfa=cmss10 scaled\magstep1
\font\sc=cmcsc10

\def\({[}
\def\){]}

 \catcode`\@=11
\def\bqed{\ifhmode\unskip\nobreak\fi\quad
  \ifmmode\blacksquare\else$\m@th\blacksquare$\fi}
\def\mex{\qopname@{mex}}

\catcode`\@=\active
\UseAMSsymbols
\vglue1.0cm
\topmatter
\title{}\endtitle
\endtopmatter

\centerline{{\tit A New Heap Game}}\bigskip\bigskip
\centerline{{\sfa Aviezri S. Fraenkel
\footnote{\tt fraenkel\@wisdom.weizmann.ac.il\qquad
http://www.wisdom.weizmann.ac.il/\~{}fraenkel}
and Dmitri Zusman\footnote{\tt dimaz\@wisdom.weizmann.ac.il}}}\bigskip
\centerline{Department of Applied Mathematics and Computer Science} 
\centerline{Weizmann Institute of Science} 
\centerline{Rehovot 76100, Israel}
\vskip1.0cm
{\bf Abstract.} Given $k\ge 3$ heaps of tokens. The moves of the $2$-player 
game introduced here are to either take a positive number of tokens from 
at most $k-1$ heaps, or to remove the {\sl same\/} positive number of 
tokens from all the $k$ heaps. We analyse this extension of Wythoff's 
game and provide a polynomial-time strategy for it. 
\vskip0.8cm
\document
\centerline{\bf 1. Introduction}\medskip

We propose the following two-player game on $k$ heaps with finitely 
many tokens, where $k\ge 3$. There are two types of moves: (i) remove 
a positive number of tokens from up to $k-1$ heaps, possibly $k-1$ 
entire heaps, or, (ii) remove the {\sl same\/} positive number of 
tokens from all the k heaps. The player making the last move wins. 

Any position in this game can be described in the following standard 
form: $(m_0,\dots,m_{k-1})$ with $0\le m_0\le\dots\le m_{k-1}$, where 
$m_i$ is the number of tokens in the $i$-th heap. 
Given any game $\Gamma$, we say informally that 
a $P$-{\it position\/} is any position $u$ of $\Gamma$ from which the 
{\it P}revious player can force a win, that is, the opponent of the player 
moving from $u$. An $N$-{\it position\/} is any position $v$ of $\Gamma$ 
from which 
the {\it N}ext player can force a win, that is, the player who moves 
from $v$. The set of all $P$-positions of $\Gamma$ is denoted by $\Cal P$, 
and the set of all $N$-positions by $\Cal N$. Denote by $F(u)$ all the 
followers of $u$, i.e., the set of all positions that can be reached 
in one move from the position $u$. It is then easy to see that:
$$\multline
\text{For every position}\ u\ \text{of}\ \Gamma\ \text{we have}\ u\in\Cal P\ 
\text{if and only if}\ F(u)\subseteq\Cal N;\\ 
\text{and}\ u\in\Cal N\ \text{if and only if}\ F(u)\cap\Cal P\ne\emptyset. 
\endmultline\tag 1$$
For $n\in\Bbb Z^0$, denote the $n$-th {\sl triangular number\/} by 
$T_n=\frac 12 n(n+1)$. We prove, 

\proclaim{\bf Theorem 1} Every\/ $P$-position of the game can be written 
in the form\/ $(T_n,m_1,\dots,m_{k-1})$, where the $(k-1)$-tuples\/ 
$(m_1,\dots,m_{k-1})$ range over all the\/ $(\!$unordered\/$)$ partitions 
of\/ $(k-1)T_n+n$ with parts of size\/ $\ge T_n$. In other words,\/ 
$\Cal P=\bigcup_{n=0}^{\infty} P_n$, where 
$$\multline
P_n=\bigl\{(T_n,m_1,\dots,m_{k-1}):\sum_{i=1}^{k-1} m_i=(k-1)T_n+n,\ \\
T_n\le m_1\le\dots\le m_{k-1},\ n\in\Bbb Z^0\bigr\}.
\endmultline\tag 2$$
\endproclaim

{\bf Example.} For $k=4$, 
$$P_n=\{(T_n,m_1,m_2,m_3):m_1+m_2+m_3=3T_n+n,\ n\in\Bbb Z^0\}.$$ 
The first few $P$-positions are: 

$P_0=\{(0,0,0,0)\}$

$P_1=\{(1,1,1,2)\}$

$P_2=\{(3,3,3,5),\ (3,3,4,4)\}$

$P_3=\{(6,6,6,9),\ (6,6,7,8),\ (6,7,7,7)\}$

$P_4=\{(10,10,10,14),\ (10,10,11,13),\ (10,10,12,12),\ (10,11,11,12)\}$

$P_5=\{(15,15,15,20),\ (15,15,16,19),\ (15,15,17,18),\ \hfill\break
\hbox{\hskip7.3cm} (15,16,16,18),\ (15,16,17,17)\}$.\medskip

\centerline{\bf 2. The Proof}\medskip
Throughout, as in (2), every $k$-tuple $(T_n,m_1,\dots,m_{k-1}),\ $ 
$(m_0,\dots,m_{k-1})$ or $(k-1)$-tuple $(m_1,\dots,m_{k-1})$ is 
arranged in nondecreasing order. Any of the first two tuples is also 
called a position (of the game) or partition (of $kT_n+n$); and the 
third is also a partition (of $(k-1)T_n+n$). The terms $m_i$ are 
called components (of the tuple) or parts (of the partition). 

\proclaim{\bf Lemma 1} Given any partition\/ $(m_1,\dots,m_{k-1})$ 
of\/ $(k-1)T_n+n$, where each part has size\/ $\ge T_n$. Then each 
part has size\/ $<T_{n+1}$.\endproclaim

\demo{\bf Proof} We have, 
$$(k-1)T_n+n-m_{k-1}=\sum_{i=1}^{k-2} m_i\ge (k-2)T_n.$$
Hence for all $i\in\{1,\dots,k-1\},\ $ $m_i\le m_{k-1}\le T_n+n=T_{n+1}-1$.
\bqed\enddemo\medskip 

\proclaim{\bf Lemma 2} Let\/ $k\ge 3$ and\/ $n\in\Bbb Z^0$. Every integer 
in the semi-closed interval\/ $t\in [T_n,T_{n+1})$ appears as a component 
in some position of\/ $P_n$. It appears in\/ $P_m$ for no\/ $m\ne n$.
\endproclaim

\demo{\bf Proof} The smallest component in $P_n$ is $T_n$, and by Lemma~1, 
the largest part cannot exceed $T_n+n=T_{n+1}-1$. Hence $t\in [T_n,T_{n+1})$ 
appears as a component in $P_m$ for no $m\ne n$. Let $t\in [T_n,T_{n+1})$, 
say $t=T_n+j,\ $ $0\le j\le n$. Then for $k\ge 3$, $T_n+j$ appears in 
the partition $\{m_1,\dots,m_{k-1}\}=\{T_n^{k-3},T_n+n-j,T_n+j\}$ of 
$(k-1)T_n+n$, where $T_n^{k-3}$ denotes $k-3$ copies of $T_n$, and so 
$T_n+j$ appears in some position of $P_n$.\bqed\enddemo\medskip 

\demo{\bf Proof of Theorem~1} It follows from (1) that it suffices to 
show two things: (I) A player moving from any position in $P_n$ lands 
in a position 
which is in $P_m$ for no $m$. (II) From any position which is in $P_m$ 
for no $m$, there is a move to some $P_n$, $n\in\Bbb Z^0$. The fact 
that (I) and (II) suffice in general for characterizing $\Cal P$ and 
$\Cal N$, is  shown in \cite{Fra$\ge$99} for the case of games without 
cycles, based on a formal definition of the $P$- and $N$-positions, 
and a proof of (1). (It is not true for cyclic games: given a digraph 
consisting of two vertices $u$ and $v$, and an edge from $u$ to $v$, 
and an edge from $v$ to $u$. Place a token on $u$. The two players 
alternate in pushing the token to a follower. The outcome is clearly 
a {\sl draw\/}, since there is no last move. However, putting 
$\Cal P=\{u\},\ $ $\Cal N=\{v\}$, satisfies (1).)

(I) Let $P_n$ be any $k$-tuple of the form (2). 
Removing tokens from up to $k-1$ heaps, including the first heap, 
results in a position $Q$ such that the first element is in $P_j$ 
for some $j<n$, yet there is a heap whose size is a component in 
$P_n$. Thus $Q\in P_m$ for no $m$ by Lemma~2. Removing tokens from 
up to $k-1$ heaps, excluding the first heap, results in a position 
$Q$ whose last $k-1$ components sum to a number $<(k-1)T_n+n$. Since, 
however, the first component is in $T_n,\ $ $Q$ is not of the form~(2). 
Hence $Q\in P_m$ for no $m$. 

So consider the move from $P_n$ which results in 
$Q=(T_n-t,m_1-t,\ldots,m_{k-1}-t)$ for some $t\in\Bbb Z^+$. If 
$Q\in P_m$ for some $m<n$, then $T_n-t=T_m$. Then 
$(T_n-t)+(m_1-t)+\dots+(m_{k-1}-t)=kT_n+n-kt=kT_m+m$. Thus, 
$0=k(T_n-T_m-t)=m-n<0$, a contradiction. Hence $Q\in P_m$ for no 
$m$. 

(II) Let $(m_0,\dots,m_{k-1})$ be any position which is in $P_m$ for no 
$m$. Since $\bigcup_{n=0}^{\infty} [T_n,T_{n+1})$ is a partition of 
$\Bbb Z^0$, we have $m_0\in [T_n,T_{n+1})$ for precisely one 
$n\in\Bbb Z^0$. Put $L=\sum_{i=1}^{k-1} m_i$. 

{\sc Case} (i).\quad $m_0=T_n$. If $L>(k-1)T_n+n$, then removing 
$L-(k-1)T_n-n$ from a suitable subset of $\{m_1,\dots,m_{k-1}\}$, 
results in a position in $P_n$. So suppose that $L<(k-1)T_n+n$. 
Then $L=(k-1)T_n+j$ for some $j\in\{0,\dots,n-1\}$. Subtracting $T_n-T_j$ 
from all components then leads to a position in $P_j$. Indeed, 
$m_0-(T_n-T_j)=T_j$, and 
$\sum_{i=1}^{k-1} \left(m_i-(T_n-T_j)\right)=(k-1)T_j+j$. 

{\sc Case} (ii).\quad $T_n<m_0<T_{n+1}$, say 
$m_0=T_n+j,\ $ $j\in\{1,\dots,n\}$. Suppose first that $L\ge(k-1)T_n+n+j$. 
If $m_1< T_{n+1}$, subtract $j$ from $m_0$ to get to $T_n$. By the first 
part of Lemma~2, $m_1$ is a part in some partition of $(k-1)T_n+n$. 
Then reduce, if necessary, a subset of the $m_i$ for $i>1$, so that 
$m_1+\sum_{i=2}^{k-1} m'_i=(k-1)T_n+n$. Here and below, $m'_i$ denotes 
$m_i$ after a suitable positive integer may have been subtracted from it. 
If $m_1\ge T_{n+1}$, then decrease $m_1$ to $T_n$. Then 
$T_n+\sum_{i\ne 1} m_i\ge T_n+j+T_n+(k-2)T_{n+1}\ge kT_n+(k-2)(n+1)+1\ge 
kT_n+n+2>kT_n+n$, since $k\ge 3$. Again by Lemma~2, $m_0$ is a part in 
some partition of $(k-1)T_n+n$. So reducing, if necessary, a subset of 
the $m_i$ for $i\ge 2$, we get $m_0+\sum_{i=2}^{k-1} m'_i=(k-1)T_n+n$. 

So consider the case $L\le (k-1)T_n+n+j$. We claim that subtracting 
$m_0-T_m$ from all components of $(m_0,\dots,m_{k-1})$ leads to a 
position in $T_m$, where $m=L-(k-1)m_0$. First note that 
$m=\sum_{i=1}^{k-1} m_i-(k-1)m_0\ge 0$, and 
$m=L-(k-1)m_0\le (k-1)T_n+n+j-(k-1)m_0=n-(k-2)j\le n-j<n$ (since $k\ge 3$), 
so $0\le m<n$, as required. Secondly, $m_0-(m_0-T_m)=T_m$, and 
$\sum_{i=1}^{k-1}\left(m_i-(m_0-T_m)\right)=L-(k-1)(m_0-T_m)=(k-1)T_m+m$. 
(Note that for $L=(k-1)T_n+n+j$ we provided two winning moves. The second 
leads to a win faster than the first.)

In conclusion, we see that $\bigcup_{i=0}^{\infty} P_i=\Cal P$.\bqed\enddemo

\centerline{\bf 3. Aspects of the Strategy}\medskip

We observe that the {\sl statement\/} of Theorem~1 tells a player 
whether or not it is possible to win by moving from any given position. 
The {\sl proof\/} of the theorem shows how to compute a winning move, 
if it exists. Together they form a {\sl strategy\/} for the game. 

The strategy can, in fact, be computed in polynomial time. 
Given any position $Q=(m_0,\dots,m_{k-1})$ of the game. Its input 
size is $\Theta\left(\sum_{i=0}^{k-1} (\log m_i)\right)$. Solving 
$m_0=n(n+1)/2$ leads to 
$n=\lfloor(\sqrt{1+8m_0}-1)/2\rfloor$. By Theorem~1, $Q\in\Cal P$ 
if and only if $m_0=T_n$, where $n=(\sqrt{1+8m_0}-1)/2$ is an integer, 
and $\sum_{i=1}^{k-1} m_i=(k-1)T_n+n$. Otherwise $Q\in\Cal N$, and 
the proof of Theorem~1 indicates how to compute a winning move to a 
$P_n$-position. All of this can be done in time which is polynomial 
in the input size. 

It is also of interest to estimate the density of the $P$-positions 
in the set of all game positions. Subtracting $T_n-1$ from each $m_i$ 
in the sum of (2), we get partitions of the form 
$$x_1+\dots +x_{k-1}=n+k-1,\quad 1\le x_1\le\dots\le x_{k-1}\le n+1,$$ 
where $x_i=m_i-(T_n-1)$. The number $p_{k-1}(n+k-1)$ of partitions 
of $n+k-1$ into $k-1$ positive integer parts is estimated in 
\cite{Hal86, \rm{Ch.~4}}. It is a polynomial of degree $k-1$ in 
$n+k-1$, whose leading term is $(n+k-1)^{k-2}/(k-2)!\ $. Thus the number 
of positions $P_n$ for $n\le N$ is estimated by 
$\pi(N)=\sum_{n=0}^N (n+k-1)^{k-2}/(k-2)!\ $. It is easy to see that 
$$\int_{-1}^N (x+k-1)^{k-2}/(k-2)!\,dx\le\pi(N)\le\int_0^{N+1} (x+k-1)^{k-2}/(k-2)!\,dx\ ,$$
leading to 
$$\frac{(N+k-1)^{k-1}-(k-2)^{k-1}}{(k-1)!}\le\pi(N)\le\frac{(N+k)^{k-1}-(k-1)^{k-1}}{(k-1)!}\ .$$

The total number of positions up to $P_N$ is the number of partitions 
of the form $m_0+\dots +m_{k-1}=n$, $0\le m_0\le\dots\le m_{k-1}$, 
where $n$ ranges from $0$ to $kT_N+N$. Adding $1$ to all the parts, 
we get partitions of the form $x_0+\dots +x_{k-1}=n+k,
\quad 1\le x_0\le\dots\le x_{k-1}\le n+k$, whose number is $p_k(n+k)$. 
As above, the total number of positions is thus 
estimated by $\nu(N)=\sum_{n=0}^{kT_N+N} (n+k)^{k-1}/(k-1)!\ $. 
Using integration as above, we get 
$$\frac{(kT_N+N+k)^k-(k-1)^k}{k!}\le\nu(N)\le\frac{(kT_N+N+k+1)^k-k^k}{k!}\ .$$
For large $N$, the ratio is thus about 
$$\frac{\pi(N)}{\nu(N)}\approx\frac{k}{kT_N+N+k}\biggl(\frac{N+k}{kT_N+N+k}\biggr)^{k-1}.$$

Dividing the numerator and denominator of the second fraction 
by $N^{k-1}$ results in $\pi(N)/\nu(N)=O(1/N^{k+1})$. We see 
that the $P$-positions are rather rare, so our game sticks to the 
majority of games in the sense of \cite{Sin81} and \cite{Sin82}. The 
rareness of $P$-positions in general, is, in fact, consistent with the 
intuition suggested by (1): a position is in $\Cal P$ if and only if 
{\sl all\/} of its followers are in $\Cal N$, whereas for a position 
to be in $\Cal N$ it suffices that one of its followers is in $\Cal P$. 
The scarcity of the $P$-positions is the reason why game strategies are 
usually specified in terms of their $P$-positions, rather than in terms 
of their $N$-positions.\medskip

\centerline{\bf 4. Epilogue}\medskip

In the heap games known to us, such as those discussed in \cite{BCG82}, 
the moves are restricted to a {\sl single\/} heap (which might, in 
special cases, be split into several subheaps). We know of two exceptions. 
One is Moore's Nim$_k$, \cite{Moo10}, where up to $k$ heaps can be 
reduced in a single move (so Nim$_1$ is ordinary Nim). The other is 
Wythoff's game, Wyt, \cite{Wyt07}, \cite{Cox53}, \cite{Dom64}, 
\cite{YaYa67}, where a move may affect up to two heaps. The motivation 
for the present note was to extend Wythoff's game to more than two heaps. 

Wyt is played on two heaps. The moves are to either remove any positive 
number of tokens from a single heap, or to remove the {\sl same\/} 
positive number of tokens from both heaps. Denoting by $(x,y)$ the 
positions of Wyt, where $x$ and $y$ denote the number of tokens in 
the two heaps with $x\le y$, the first eleven $P$-positions are listed 
in Table~1. The reader may wish to guess the next few entries of the table 
before reading on. 

\midinsert
\topcaption{Table 1. {\rm The first few $P$-positions of Wyt.}}\endcaption
$$\vcenter{\lineskip 0pt\halign{\hfill #\qquad\qquad&\hfill #\qquad
\qquad&\hfill #\cr
\noalign{\hrule\vskip 5pt}
$n$&$A_n$&$B_n$\cr
\noalign{\vskip 5pt\hrule\vskip 5pt}
0&0&0\cr 1&1&2\cr 2&3&5\cr 3&4&7\cr 4&6&10\cr 5&8&13\cr 6&9&15\cr
7&11&18\cr 8&12&20\cr 9&14&23\cr 10&16&26\cr
\noalign{\vskip 5pt\hrule}}}$$
\endinsert\par
For any finite subset $S\subset\Bbb Z^0$, define the {\it M}inimum 
{\it EX}cluded value of $S$ as follows: $\mex S=\min\Bbb Z^0\setminus S$
= least nonnegative integer not in $S$ \cite{BCG82}. Note that if 
$S=\emptyset$, then $\mex S=0$. The general structure of Table~1 is given by: 
$$A_n= \mex\{A_i, B_i: 0\le i<n\},\ B_n= A_n+n\quad (n\in\Bbb Z^0).$$
 
Since the input size of Wyt is {\sl succinct\/}, namely $\Theta (\log (x+y))$, 
one can see that the above 
characterization of the $P$-positions implies a strategy which is 
exponential. A polynomial strategy for Wyt can be based on the observation 
that $A_n=\lfloor n\alpha\rfloor$, $B_n=\lfloor n\beta\rfloor$, where 
$\alpha=(1+\sqrt{5})/2$ is the golden section, $\beta=(3+\sqrt{5})/2$. 
Another polynomial strategy depends on a special numeration system 
whose basis elements are the numerators of the simple continued fraction 
expansion of $\alpha$. These three strategies can be generalized to 
Wyt$_a$, proposed and analysed in \cite{Fra82}, where $a\in\Bbb Z^+$ 
is a parameter of the game. The moves are as in Wyt, except that the 
second type of move is to remove say $k>0$ and $l>0$ from the two heaps 
subject to $|k-l|<a$. Clearly Wyt$_1$ is Wyt. 

The generalization of Wyt to more than two heaps was a long sought-after 
problem. In \cite{Fra96} it is shown that the natural generalization 
to the case of $k\ge 2$ heaps is to either remove any positive number 
of tokens from a single heap, or say $l_1,\dots,l_k$ from all of 
them simultaneously, where the $l_i$ are nonnegative integers 
with $\sum_{i=1}^k l_i>0$ and $l_1\oplus\dots\oplus l_k=0$, and where 
$\oplus$ denotes Nim-sum (also known as addition over GF(2), or XOR). 
In particular, the case $k=2$ is Wyt. But the actual computation 
of the $P$-values seems to be difficult. 

The heap-game considered here is a generalization of the {\sl moves\/} 
of Wyt, but not of its strategy. In fact, it doesn't specialize to the 
case $k=2$; we used the fact that $k\ge 3$ in several places of the proof. 
However, the $P$-positions of the present game have a compact form, 
the exhibition of which was the purpose of this note. 

We remark finally that the {\sl Sprague-Grundy\/} function $g$ of a game 
provides a strategy for the {\sl sum\/} of several games. The computation 
of $g$ for Nim$_k$, $k\ge2$, and Wyt$_a$, $a\ge 1$ seems to be difficult. 
It would be of interest to compute the $g$-function for the present 
game. Perhaps this is also difficult.\medskip

{\bf Acknowledgment.} We thank Uriel Feige for a simplification in the 
statement and proof of Theorem~1.\medskip

\centerline{\bf References}\medskip

1. [BCG82] E. R. Berlekamp, J. H. Conway and R. K. Guy \(1982\), {\it Winning
Ways\/} (two volumes), Academic Press, London.

2. [Cox53] H. S. M. Coxeter \(1953\), The golden section, phyllotaxis 
and Wythoff's game, {\it Scripta Math.\/} {\bf 19}, 135--143.

3. [Dom64] A. P. Domoryad \(1964\), {\it Mathematical Games and Pastimes\/}
(translated by H. Moss), Pergamon Press, Oxford.

4. [Fra82] A. S. Fraenkel \(1982\), How to beat your Wythoff games' 
opponent on three fronts, {\it Amer. Math. Monthly\/} {\bf 89}, 353--361.

5. [Fra96] A. S. Fraenkel \(1996\), Scenic trails ascending from sea-level Nim 
to alpine chess, in: {\it Games of No Chance}, Proc. MSRI Workshop 
on Combinatorial Games, July, 1994, Berkeley, CA (R. J. Nowakowski, 
ed.), MSRI Publ. Vol.~29, Cambridge University Press, Cambridge, pp.~13--42. 

6. [Fra$\ge$99] A. S. Fraenkel [$\ge$1999], {\it Adventures in Games and 
Computational Complexity}, to appear in Graduate Studies in Mathematics,
Amer. Math. Soc., Providence, RI.

7. [Hal86] M. Hall \(1986\), {\it Combinatorial Theory\/}, $2$nd edition, 
Wiley, New York. 

8. [Moo10] E. H. Moore \(1909--1910\), A generalization of the game called 
nim, {\it Ann. of Math.} {\bf 11} (Ser.~2), 93--94.

9. [Sin81] D.\ Singmaster \(1981\), Almost all games are first person 
games, {\it Eureka\/} {\bf 41}, 33--37.

10. [Sin82] D.\ Singmaster \(1982\), Almost all partizan games are first 
person and almost all impartial games are maximal, {\it J.\ Combin.\ Inform.\ 
System Sci.\/} {\bf 7}, 270--274.

11. [Wyt07] W. A. Wythoff \(1907\), A modification of the game of Nim, 
{\it Nieuw Arch. Wisk.\/} {\bf 7}, 199--202.

12. [YaYa67] A. M. Yaglom and I. M. Yaglom \(1967\), {\it Challenging 
Mathematical Problems with Elementary Solutions\/}, translated by 
J. McCawley, Jr., revised and edited by B. Gordon, Vol.~II, Holden-Day, 
San Francisco.

\enddocument